\documentclass[12pt,a4paper]{article}
\usepackage[latin1]{inputenc}
\usepackage{amssymb}
\usepackage{amsmath}
\usepackage{stmaryrd}
\usepackage[english]{babel}
\usepackage{cases}
\usepackage{fancyhdr}
\usepackage{graphicx}
\usepackage{color}
\usepackage{setspace}

\usepackage[latin1]{inputenc}
\usepackage{color}
\usepackage{latexsym}
\usepackage{amssymb}
\usepackage{graphicx}
\usepackage{enumerate}
\usepackage{extarrows}
\newtheorem{Theorem}{Theorem}[part]

\newtheorem{Remark}{Remark}[part]

\def \2{\vspace{2mm}}

\def \nn {\nonumber}

\def \Sum{\displaystyle\sum}

\def \Frac{\displaystyle\frac}

\def \Lim{\displaystyle\lim}

\def \R{\mathbb{R}}

\def \Z{\mathbb{Z}}
\def \E{\mathbb{E}}

\def \Ec{{\cal E}}

 \def \Nc{{\cal N}}

\def \Hf    {\mathfrak{H}}

\def \o {\otimes }

\def \ep{\hbox{ }\hfill$\Box$}

\def\be{\begin{eqnarray}}
\def\ee{\end{eqnarray}}

\def\b*{\begin{eqnarray*}}
\def\e*{\end{eqnarray*}}
\def \nn{\nonumber }

\newcommand{\nodate}{\date{}}\nodate

\begin{document}

\title{A Simple  Proof of Berry-Esséen Bounds for the Quadratic Variation of the Subfractional Brownian Motion }
\maketitle
 \begin{center}
Soufiane Aazizi \footnote{Department of Mathematics, Faculty of Sciences Semlalia Cadi Ayyad University, B.P. 2390 Marrakesh, Morocco. Email: {\tt
aazizi.soufiane@gmail.com}

$^*$The author is supported by the Marie Curie Initial Training Network (ITN) project: ``Deterministic and Stochastic Controlled Systems and Application",
FP7-PEOPLE-2007-1-1-ITN, No. 213841-2.}
\\
{\it  Université Cadi Ayyad}
\end{center}

\begin{abstract}
We give a simple technic  to derive the Berry-Esséen bounds for the quadratic variation of the subfractional Brownian motion (subfBm). Our approach has two main ingredients: ($i$) bounding from above the covariance of quadratic variation of subfBm by the covariance of the quadratic variation of fractional Brownian motion (fBm); and ($ii$) using the existing results on fBm in \cite{BN08,NP09,N12}. As a result, we obtain simple and direct proof to derive the rate of convergence of quadratic variation of subfBm. In addition, we also improve this rate of convergence to meet the one of fractional Brownian motion in \cite{N12}.

\end{abstract}

\vspace{9mm}

\noindent {\bf Key words~:} Fractional Brownian motion, Malliavin calculus, Kolmogorov distance, Subfractional Brownian motion, Stein method, Quadratic variation.

\setcounter{equation}{0} \setcounter{Assumption}{0}
\setcounter{Theorem}{0} \setcounter{Proposition}{0}
\setcounter{Corollary}{0} \setcounter{Lemma}{0}
\setcounter{Definition}{0} \setcounter{Remark}{0}

\section{Introduction and preliminaries }
The following result, proved in \cite{T11}, shows the convergence of quadratic variation of subfractional Brownian motion (subfBm in short) to a centered reduced normal variable, the author also provides its rate of convergence.
Let $S = ( S_t, \, t\geq 0 )$ be a  subfractional Brownian motion, and define
 \b*
 Z_n &=& \sum_{k=0 }^{n-1}n^{2H}\left[ (S_{(k+1)/n}-S_{k/n})^2-Var\left(S_{(k+1)/n}-S_{k/n}\right)\right], \quad n\geq1.
 \e*

\begin{Theorem}{(\textbf{Tudor 2011})}
Let $N$ be a standard Gaussian random variable ($N \sim N(0, 1)$) and suppose that $H \in (0, \frac34]$. Then $\frac{Z_n}{Var(Z_n)}$ converges in distribution to $N$ and the following Berry-Esséen bounds hold for every $n \ge 1$,
\b*
    d_{Kol}\left(\frac{Z_n}{Var(Z_n)},N\right) &\le& c_{H} \times
    \left\{
     \begin{array}{ll}
       n^{-\frac12} , & H\in \left(0, \frac12\right), \\
       ~~\\
       n^{2H-\frac32}, & H \in \left[\frac12, \frac34\right),  \\
       ~~\\
        \frac{1}{\sqrt{\log n}}, & H  =\frac{3}{4},
 \end{array}
   \right.
\e*
where $c_H$ is a constant depending only on $H$.
\end{Theorem}

In \cite{T11}, the proof uses stein method and malliavin calculus, based on the idea developed in \cite{BN08,NP09} for the case of fractional Brownian motion (fBm in short), which leads to the same rate of convergence.  Recently, \cite{N12} used the convolution product of two sequences which improve clearly the rate of convergence of the fBm. The natural question imposes itself, it is possible to obtain a rate of convergence of subfBm similar to the one proved by \cite{N12} for the fBm? \\

The goal of this paper, is to improve the rate of convergence of the subfBm so that we have at least the same one as the fBm. To perform our calculation, we will mainly follow the idea taken from \cite{N12}. With the proof of \cite{NP09} and \cite{N12} in hand, we will show  how we can retrieve the result of \cite{T11}, and how we can improve this result to reach the one of fBm in \cite{N12}. We claim the main result of this paper:
\begin{Theorem} \label{CLT}
Let $N \sim \Nc(0,1)$, there exist a constant $c_H$ depending only on $H$, such that for every $n \ge 1$,
 \b*
    d_{Kol}\left(\frac{Z_n}{Var(Z_n)},N\right) &\le& c_{H} \times
    \left\{
     \begin{array}{ll}
       \frac{1}{\sqrt{n}} , & H\in \left(0, \frac58\right), \\
       ~~\\
       \frac{(\log n)^{3/2}}{\sqrt{n}}, & H=\frac58, \\
       ~~\\
       n^{4H-3}, & H  \in \left(\frac58 ,\frac{3}{4}\right),\\
       ~~\\
       \frac{1}{\log n}, & H  =\frac{3}{4}.
 \end{array}
   \right.
\e*
\end{Theorem}
The subfBm $S = ( S_t, \, t\geq 0 )$ with parameters $H\in (0, 1)$, is defined on some probability space $(\Omega,\mathcal{F},P)$ (Here, and everywhere else, we do assume that $\mathcal{F}$ is the sigma-field generated by $S$). This means that $S$ is a centered Gaussian process with  covariance
 \be\label{cov}
\E[S_sS_t]&=&R_H(s,t)=s^{2H}+ t^{2H} -\frac12\left[(s+t)^{2H} + |t-s|^{2H} \right], \quad s, t \ge 0.
\ee
We recall briefly some important tools of Malliavin calculus used throughout this paper.  We mean by  $\Hf $ a real separable Hilbert space defined as follows: (i) denote by $\Ec$ the set of all $\R$-valued  functions on $[0,\infty)$, (ii) define $\Hf$ as the Hilbert space obtained by closing $\Ec$
with respect to the scalar product
 \b*
 \langle1_{[0, s]}, 1_{[0, t]}\rangle_\Hf&=&R_H(s,t).
 \e*
 For every $q\geq 1$, let $\mathcal{H}_q$ be the $q^{\mbox{th}}$ Wiener chaos of $X$, that is, the closed linear subspace of $L^{2}(\Omega)$ generated by the random variables
$\{H_{q}\left( X\left( h\right) \right) ,h\in \Hf,\| h\| _{\Hf}=1\}$, where $H_{q}$ is the $q^{\mbox{th}}$ Hermite polynomial defined as $H_{q}(x)=(-1)^q \;
e^{\frac{x^2}{2}}\;\Frac{d^q}{dx^q}(e^{-\frac{x^2}{2}})$. The mapping $I_{q}(h^{\otimes q})=H_{q}\left( X\left( h\right) \right) $ provides a linear isometry between
the symmetric tensor product $\Hf^{\odot q}$ (equipped with the modified norm $\|\cdot\|_{\Hf^{\odot q}}= \sqrt{q!}  \;\|\cdot\|_{\Hf^{\otimes q}}$) and
$\mathcal{H}_q$.
Specifically, for all $f,g\in\Hf^{\odot q}$ and $q\geq 1$, one has
\be\label{isoint}
\E\big[I_q(f)I_q(g)\big]&=&q! \;\langle f,g\rangle_{\mathcal{H}^{\otimes q}}.
\ee
Let $\{e_k, k\geq 1\}$ be a complete orthonormal system in $\Hf$. Given $f\in\Hf^{\odot p}$ and $g\in\Hf^{\odot q}$, for every $r=0, \dots, p\wedge q$, the $r^{\mbox{th}}$ contraction of $f$ and $g$ is the element of $\Hf^{\otimes(p+q-2r)}$ defined as
$$
f\otimes_r g=\sum_{i_1=1, \dots, i_r=1}^{\infty}
\langle f,e_{i_1}\otimes\cdots\otimes e_{i_r}\rangle_{\Hf^{\otimes r}}
\otimes
\langle g,e_{i_1}\otimes\cdots\otimes e_{i_r}\rangle_{\Hf^{\otimes r}}.
$$
In particular, note that $f\otimes_0g=f\otimes g$ and when $p=q$,
that $f\otimes_pg=\langle f, g\rangle_{\Hf^{\otimes p}}$.  Since, in general, the contraction $f\otimes_rg$ is not necessarily symmetric,
we denote its symmetrization by $f\widetilde\otimes_rg \in \Hf^{\odot(p+q-2r)}$. The following formula is useful to compute the product of such multiple
integrals: if $f\in \Hf^{\odot p}$ and $g\in \Hf^{\odot q}$, then
\be
\label{eq:multiplication}
I_p(f)\;I_q(g)&=&\sum_{r=0}^{p\wedge q} r! \left(\!\!\begin{array}{c}p\\r\end{array}\!\!\right) \left(\!\!\begin{array}{c}q\\r\end{array}\!\!\right)
I_{p+q-2r}(f\widetilde\otimes_rg).
\ee
We will use the notation $ \delta_{k/n}=1_{[k/n,(k+1)/n]}$, and we send the reader to \cite{N06} for more details on Malliavin calculus.\\
Now, by  self-similarity property of $S$ and (\ref{cov}) we deduce  for $k \le l$
 \b*
 n^{2H}\langle\delta_{k/n},\delta_{l/n}\rangle_\Hf&=& n^{2H}\E\left(\left(S_{(k+1)/n}-S_{k/n}\right)\left(S_{(l+1)/n}-S_{l/n}\right)\right)\\
 &=&\E\left( (S_{k+1}-S_k)(S_{l+1}-S_l)\right)\\
 &=&(k+l+1)^{2H} -\frac12(k+l+2)^{2H}-\frac12(k+l)^{2H} \\
 && - (l-k)^{2H}  +\frac12(l-1-k)^{2H}+\frac12(l+1-k)^{2H}\\
 &=&\frac12\rho(l-k)-\frac12\rho(l+k+1),
 \e*
where $\rho(r)=|r+1|^{2H}+|r-1|^{2H}-2|r|^{2H},\quad r\in \mathbb{Z }$.
\vspace{2mm}\\
So that, we have the relation
\be
\left| n^{2H}\langle\delta_{k/n},\delta_{l/n}\rangle_\Hf\right|&=& \frac12|\rho(l-k)-\rho(l+k+1)| \nn\\
&\le& |\rho(l-k)| \label{Relation Cov rho},
\ee
since the function $r \rightarrow |\rho(r)|$ is nonincreasing. In fact, we can write $\rho$ as
\b*
\rho(r)&=&f(r+1)-f(r),
\e*
where
$f(r):=|r+1|^{2H}-|r|^{2H}$. It follows that:\\
For $H \ge \frac12: f' >0 \rightarrow \rho >0$ and $f'' < 0 \rightarrow \rho  \searrow$, which implies that $|\rho|$ is nonincreasing.\vspace{1.5mm}\\
For $H \le \frac12: f' <0 \rightarrow \rho < 0$ and $f'' > 0 \rightarrow \rho  \nearrow$, which implies that $|\rho|$ is nonincreasing.
\vspace{2mm}
\\
With inequality (\ref{Relation Cov rho}) in hand, it is now straightforward to obtain Theorem \ref{CLT}. Hence, we can write the quadratic variation of $S$, with respect to a subdivision $\mathcal{\pi}_n = \{0< \frac{1}{n} <\frac{2}{n}<\ldots < 1\}$ of $[0, 1]$, as
follows
 \be
Z_n&=&\Sum_{k=0}^{n-1}\left[ n^{2H}\left(S_{(k+1)/n}-S_{k/n}\right)^2 -1+\frac12\rho(2k+1)\right]\nn\\
&=&\Sum_{k=0}^{n-1}\left[ n^{2H}\left(I_1(\delta_{k/n})\right)^2 -1+\frac12\rho(2k+1)\right]\nn\\
&=&I_2\left(\underbrace{\Sum_{k=0}^{n-1}n^{2H}\delta^{\o 2}_{k/n}}_{g_n}\right). \label{Z_n}
\ee
Thus, we can write the correct
renormalization of $Z_n$  as follows,
\be \label{V_n second
expression}V_n=\frac{Z_n}{\sqrt{Var(Z_n)}}=\frac{I_2(g_n)}{\sqrt{Var(Z_n)}}.
\ee
\section{Proof of Theorem \ref{CLT}}
In the first step, we show that $\frac{Var(Z_n)}{n}$  and $\frac{Var(Z_n)}{n\log n}$  have a limit. Therefore, we have
\b*
\frac{Var(Z_n)}{n}&=&n^{-1}\E[I_2^2(g_n)]=2\|g_n\|^2_{\Hf^{\o 2}}\\
    &=&2n^{4H-1}\Sum_{k,l=0}^{n-1} \langle \delta^{\o 2}_{k/n},\delta^{\o 2}_{l/n}\rangle_{\Hf^{\o 2}}
    =2n^{4H-1}\Sum_{k,l=0}^{n-1} \langle \delta_{k/n},\delta_{l/n}\rangle_{\Hf}^2 \\
    &=&\frac{1}{2n}\Sum_{k,l=0}^{n-1} \left| \rho(l-k)-\rho(l+k+1)\right|^2\\
    &=&\frac{1}{2n}\Sum_{k,l=0}^{n-1} \rho^2(l-k)+ \frac{1}{2n}\Sum_{k,l=0}^{n-1} \rho^2(l+k+1)-\frac1n\Sum_{k,l=0}^{n-1} \rho(l-k)\;\rho(l+k+1).
\e*
As in the proof of \cite[Theorem 5.6]{N12}, we have for $H < \frac34$
\be
\Lim_{n\rightarrow \infty }\frac1n\Sum_{k,l=0}^{n-1} \rho^2(l-k) &=& \Sum_{r\in \Z}\rho^2(r) \label{rho(l-k)}.
\ee
On the other hand
\b*
\Sum_{k,l=0}^{n-1} \rho^2(l+k+1)&=& \sum_{|r|<2n-1}(r+1)\;\rho^2(r+1).
\e*
Assume that $H<\frac34$ and write
\b*
\frac1n\Sum_{k,l=0}^{n-1} \rho^2(l+k+1)&=& \sum_{r\in \Z}\rho^2(r+1)\;\frac{r+1}{n} \; \mathbf{1}_{\{| r| < 2n-1\}}.
\e*
From \cite [Lemma 4.3]{NP09}, we have for any $\alpha \in \R$ we have
\be
\Sum_{k=1}^{n-1} k^\alpha &\trianglelefteqslant& 1+ n^{\alpha+1},
\ee
where the notation $a_n \trianglelefteqslant b_n$ means that $\sup_{n\ge1} |a_n|/|b_n| < \infty$. Combined with the fact that the function $\rho$ behaves asymptotically as
\b*\rho(r)=2HK(2HK-1)\;|r|^{2HK-2},\quad |r|\rightarrow\infty.
\e*
Leads to
\be
\lim_{n \rightarrow \infty }\frac1n\Sum_{k,l=0}^{n-1} \rho^2(l+k+1)&=& 0. \label{rho(l+k+1)}
\ee
Finally, by (\ref{rho(l-k)}) and (\ref{rho(l+k+1)}), together with Cauchy Schwartz inequality
\be
\frac1n\Sum_{k,l=0}^{n-1} |\rho(l-k)|~ |\rho(l+k+1)|
&\le& \left(\frac1n\Sum_{k,l=0}^{n-1} \rho^2(l-k) \right)^\frac12\left( \frac1n\Sum_{k,l=0}^{n-1} \rho^2(l+k+1)\right)^\frac12 \nn\\
&\longrightarrow& 0, \quad \mbox{ as } n\rightarrow 0 \label{rho(l-k)rho(l+k+1)}.
\ee
Combining (\ref{rho(l-k)}), (\ref{rho(l+k+1)}) and (\ref{rho(l-k)rho(l+k+1)}) we conclude that
\be
\Lim_{n\rightarrow \infty }\frac{Var(Z_n)}{n}&=& \frac{1}{2}\Sum_{r\in \Z}\rho^2(r) \label{VarZn/n}.
\ee
Assume now $H=\frac34$. Following similar argument as above we have
\b*
\frac{Var(Z_n)}{n\log(n)}    &=&\frac{1}{n\log(n)}\;\Sum_{k,l=0}^{n-1} \rho^2(l-k)+ \frac{1}{n\log(n)}\;\Sum_{k,l=0}^{n-1} \rho^2(l+k+1)\\
&&-\frac{2}{n\log(n)}\;\Sum_{k,l=0}^{n-1} \rho(l-k)\;\rho(l+k+1).
\e*
Again from the proof of \cite[Theorem 5.6]{N12}, we have
\be
\lim_{n\rightarrow \infty }\;\frac{1}{n\log(n)}\;\Sum_{k,l=0}^{n-1} \rho^2(l-k) &=& \frac{9}{32}. \label{rho(l-k)/ nlog(n)}
\ee
From other side, we have $\rho^2(r) \sim \frac{9}{64|r|}$ as $|r| \rightarrow \infty$. Implying in turn
\b*
\Sum_{k,l=0}^{n-1} \rho^2(l+k+1)&=&\sum_{r=1}^{2n-1}r\rho^2(r)\sim \frac{9}{64}\sum_{r=1}^{2n-1}1\sim (2n-1)\frac{9}{64},\quad \mbox{ as } n\rightarrow \infty.
\e*
Hence, we have
\be
\lim_{n\rightarrow \infty }\;\frac{1}{n\log(n)}\;\Sum_{k,l=0}^{n-1} \rho^2(l+k+1)&=&0. \label{rho(l+k+1)/ nlog(n)}
\ee
Similarly to (\ref{rho(l-k)rho(l+k+1)}), we obtain by (\ref{rho(l-k)/ nlog(n)}), (\ref{rho(l+k+1)/ nlog(n)}) and Cauchy Schwartz
\be
\lim_{n\rightarrow \infty } \;\frac{1}{n\log(n)}\;\Sum_{k,l=0}^{n-1} \rho(l-k)\;\rho(l+k+1)&=&0. \label{rho(l-k)rho(l+k+1)/nlogn}
\ee
Combining (\ref{rho(l-k)/ nlog(n)}), (\ref{rho(l+k+1)/ nlog(n)}) and (\ref{rho(l-k)rho(l+k+1)/nlogn}) we deduce that
\be
\frac{Var(Z_n)}{n\log(n)}    &=&\frac{9}{64} . \label{VarZn/n logn}
\ee
Let us now derive the explicit bounds. From (\ref{Z_n}), multiplication formula (\ref{eq:multiplication}) and the fact that
$\E \|DZ_n\|^2_\Hf=2Var(Z_n)$, we obtain
\b*
\frac{1}{2}\|DV_n\|^2_\Hf -1 &=&\frac{2n^{4H}}{Var(Z_n)}\Sum_{k,l=0}^{n-1}I_2(\delta_{k/n}\widetilde{\o}\delta_{l/n})\;\langle \delta_{k/n},\delta_{l/n}\rangle_\Hf.
\e*
It follows  by (\ref{Relation Cov rho}) that
\be &\E&\left[\left(\frac{1}{2}\left\|DV_n\right\|^2_\Hf - 1\right)^2\right] \nn\\
&=& \frac{4n^{8HK}}{Var^2(Z_n)} \E\left[\left( \Sum_{k,l=0}^{n-1}I_2(\delta_{k/n}\widetilde{\o }\delta_{l/n})\;\langle \delta_{k/n},\delta_{l/n}\rangle_\Hf\right)^2\right]\nn\\
&=&\frac{8n^{8HK}}{Var^2(Z_n)}\Sum_{i,j,k,l=0}^{n-1}\langle \delta_{k/n},\delta_{l/n}\rangle_\Hf \;\langle \delta_{k/n},\delta_{l/n}\rangle_\Hf \; \langle
\delta_{k/n}\widetilde{\o}\delta_{l/n},\delta_{k/n}\widetilde{\o}\delta_{l/n}\rangle_{\Hf^{\o 2}} \nn \\
&=&\frac{4n^{8HK}}{Var^2(Z_n)}\Sum_{i,j,k,l=0}^{n-1}\langle \delta_{i/n},\delta_{j/n}\rangle_\Hf \; \langle \delta_{k/n},\delta_{l/n}\rangle_\Hf \;\Big(\langle \delta_{i/n} ,\delta_{k/n}\rangle_{\Hf} \; \langle \delta_{j/n} ,\delta_{l/n}\rangle_{\Hf}\nn\\
&&\hspace{75mm}+\langle \delta_{i/n} ,\delta_{l/n}\rangle_{\Hf}\langle
\delta_{j/n} ,\delta_{k/n}\rangle_{\Hf}\Big)\nn\\
&=&\frac{8n^{8HK}}{Var^2(Z_n)}\ \Sum_{i,j,k,l=0}^{n-1}\; \langle \delta_{i/n},\delta_{j/n}\rangle_\Hf
  \; \langle \delta_{i/n},\delta_{k/n}\rangle_\Hf \;\langle \delta_{k/n} ,\delta_{l/n}\rangle_{\Hf} \;\langle \delta_{j/n} ,\delta_{l/n}\rangle_{\Hf} \nn \\
&\le&\frac{8n^2}{Var^2(Z_n)}\,\,\frac{1}{n^2}\Sum_{i,j,k,l=0}^{n-1} |\rho(i-j)|\;|\rho(i-k)|\;|\rho(k-l)|\;|\rho(j-l)|\label{DV_n and A(n)}.
\ee
Then, combining the convergence (\ref{VarZn/n}) and (\ref{VarZn/n logn}) together with inequality  (\ref{DV_n and A(n)}),  the rest of the proof is now similar to the one of Theorem 5.6 in \cite{N12}.
\ep

\begin{Remark}
To retrieve the result of Tudor \cite{T11}, we start from equality (\ref{VarZn/n logn}) and we follow the same steps as in the proof of Theorem 4.1 in \cite{NP09}.
\end{Remark}

\end{document}